\input amstex
\input amsppt.sty   
\hsize 30pc
\vsize 47pc
\def\nmb#1#2{#2}         
\def\cit#1#2{\ifx#1!\cite{#2}\else#2\fi}
\def\idx{}               
\def\ign#1{}             

\define\rh{\rho}

\define\ph{\varphi}

\define\ps{\psi}
\define\om{\omega}

\define\La{\Lambda}
\define\Si{\Sigma}

\define\Om{\Omega}

\define\row#1#2#3{#1_{#2},\ldots,#1_{#3}}
\define\x{\times}

\def\today{\ifcase\month\or
 January\or February\or March\or April\or May\or June\or
 July\or August\or September\or October\or November\or December\fi
 \space\number\day, \number\year}
\topmatter
\title  
Basic differential forms for actions of Lie groups, II
\endtitle
\author 
Peter W\. Michor  \endauthor
\affil
Erwin Schr\"odinger International Institute of Mathematical Physics, 
Wien, Austria \\
Institut f\"ur Mathematik, Universit\"at Wien, Austria
\endaffil
\address
P\. W\. Michor: Institut f\"ur Mathematik, Universit\"at Wien,
Strudlhofgasse 4, A-1090 Wien, Austria
\endaddress
\email Peter.Michor\@univie.ac.at \endemail
\thanks Supported by Project P 10037--PHY
of `Fonds zur F\"orderung der wissenschaftlichen Forschung'.
\endthanks
\subjclass Orbits, sections, basic differential forms\endsubjclass
\abstract 
The assumption in the main result of \cit!{2} is removed
\endabstract
\endtopmatter


\document
Let $G$ be a Lie group which acts isometrically on a Riemannian 
manifold $M$. A section of the Riemannian $G$-manifold $M$ is a 
closed submanifold $\Si$ which meets each orbit orthogonally. In this 
situation the trace on $\Si$ of the $G$-action is a discrete group 
action by the generalized Weyl group $W(\Si)=N_G(\Si)/Z_G(\Si)$, 
where $N_G(\Si):= \{g\in G: g.\Si=\Si\}$ and $Z_G(\Si):= 
\{g\in G: g.s=s \text{ for all }s\in \Si\}$.
A differential form $\ph\in \Om^p(M)$ is called 
\idx{\it $G$-invariant} if $g^*\ph=\ph$ for all $g\in G$ and 
\idx{\it horizontal} if $\ph$ kills each vector tangent to a 
$G$-orbit.
We denote by $\Om^p_{\text{hor}}(M)^G$ the space of all horizontal 
$G$-invariant $p$-forms on $M$ which are also called \idx{\it basic 
forms}.

In the paper \cit!{2} it was shown that for a 
proper isometric action of a Lie group $G$ on a smooth Riemannian 
manifold $M$ admitting a section $\Si$
the restriction of differential forms induces an isomorphism
$$\Om^p_{\text{hor}}(M)^G @>{\cong}>> \Om^p(\Si)^{W(\Si)}$$
between the space of horizontal $G$-invariant differential forms on 
$M$ and the space of all differential forms on $\Si$ which are 
invariant under the action of the generalized Weyl group $W(\Si)$ of 
the section $\Si$, under the following assumption:
\block
     For each $x\in\Si$ the slice representation 
     $G_x\to O(T_x(G.x)^\bot)$ has a generalized Weyl group which 
     is a reflection group.
\endblock
In this paper we will show that this result holds in general, without 
any assumption. Notation is as in \cit!{2}, which is used throughout.
For more information on $G$-manifolds with sections see the seminal 
paper \cit!{3}.

\subhead\nmb0{1}. Polar representations \endsubhead
Let $G$ be a compact Lie group and let $\rh:G\to GL(V)$ be an 
orthogonal representation in a finite dimensional real vector space 
$V$ which admits a section $\Si$. Then the section turns out to be a 
linear subspace and the representation is called a \idx{\it polar 
representation}, following Dadok \cit!{1}, who gave a complete 
classification of all polar representations of connected Lie groups. 

\proclaim{Theorem}
Let $\rh:G\to O(V)$ be a polar orthogonal representation of 
a compact Lie group $G$, with section $\Si$ and generalized Weyl 
group $W=W(\Si)$. 
Let $B\subset V$ be an open ball centered at 0.

Then the restriction of differential forms induces an isomorphism
$$\Om^p_{\text{hor}}(B)^G @>{\cong}>> \Om^p(\Si\cap B)^{W(\Si)}.$$
\endproclaim

\demo{Proof} 
We only treat the case $B=V$. The restriction to an open ball can 
be proved as in \cit!{2}, 3.8.
Let $i:\Si\to V$ be the embedding. 
It is easy to see (and proved in \cit!{2}, 2.4) that the restriction 
$i^*:\Om^p_{\text{hor}}(V)^G \to \Om^p(\Si)^{W(G)}$ 
is injective, so it remains to prove surjectivity.
Let $G_0$ be the connected component of $G$. From \cit!{1}, lemma 1 
one concludes:
\block
{\sl A subspace $\Si$ of $V$ is a section for $G$ if and only if it 
is a section for $G_0$. Thus $\rh$ is a polar representation for $G$ 
if and only if it is a polar representation for $G_0$.}
\endblock
The generalized Weyl groups of $\Si$ with respect to $G$ and to $G_0$ 
are related by 
$$
W(G_0)=N_{G_0}(\Si)/Z_{G_0}(\Si)\subset W(G)=N_G(\Si)/Z_G(\Si), 
$$
since $Z_G(\Si)\cap N_{G_0}(\Si)=Z_{G_0}(\Si)$.

Let $\om\in \Om^p(\Si)^{W(G)}\subset \Om^p(\Si)^{W(G_0)}$. Since 
$G_0$ is connected the generalized Weyl group $W(G_0)$ is generated 
by reflections (a Coxeter group) by \cit!{1}, remark after 
proposition 6. Thus we may 
apply \cit!{2}, theorem 3.7, which asserts that then 
$$i^*:\Om^p_{\text{hor}}(V)^{G_0} @>{\cong}>> \Om^p(\Si)^{W(G_0)}$$
is an isomorphism, and we get $\ph\in \Om^p_{\text{hor}}(V)^{G_0}$ 
with $i^*\ph=\om$. Let us consider 
$$\ps := \int_G g^*\ph\,dg\in \Om^p_{\text{hor}}(V)^G,$$
where $dg$ denotes Haar measure on $G$.
In order to show that $i^*\ps=\om$ it suffices to check that 
$i^*g^*\ph=\om$ for each $g\in G$.
Now $g(\Si)$ is again a section of $G$, thus also of $G_0$. Since any 
two sections are related by an element of the group, there exists 
$h\in G_0$ such that $hg(\Si)=\Si$. Then $hg\in N_{G}(\Si)$ and we 
denote by $[hg]$ the coset in $W(G)$, and we may compute as follows:
$$\align
(i^*g^*\ph)_x &=  (g^*\ph)_x.\La^pTi = \ph_{g(x)}.\La^pTg.\La^pTi\\
&= (h^*\ph)_{g(x)}.\La^pTg.\La^pTi,
     \quad\text{ since }\ph\in\Om^p_{\text{hor}}(M)^{G_0}\\
&= \ph_{hg(x)}.\La^pT(hg).\La^pTi
     = \ph_{i[hg](x)}.\La^pTi.\La^pT([hg])\\
&= (i^*\ph)_{[hg](x)}.\La^pT([hg]) \\
&= \om_{[hg](x)}.\La^pT([hg])= [hg]^*\om =\om.\qed
\endalign$$
\enddemo

\proclaim{\nmb0{2}. Theorem} Let $M\x G\to M$ be a proper 
isometric right action of a Lie group $G$ on a smooth Riemannian 
manifold $M$, which admits a section $\Si$. 

Then the restriction of 
differential forms induces an isomorphism
$$\Om^p_{\text{hor}}(M)^G @>{\cong}>> \Om^p(\Si)^{W(\Si)}$$
between the space of horizontal $G$-invariant differential forms on 
$M$ and the space of all differential forms on $\Si$ which are 
invariant under the action of the generalized Weyl group $W(\Si)$ of 
the section $\Si$.
\endproclaim

This is the Main Theorem 2.4 of \cit!{2}, without the assumption made 
there.

\demo{Proof}
Injectivity is proved in \cit!{2}, 2.4, without using the assumption. 
Surjectivity can be proved as in \cit!{2}, section 4, where one 
replaces the use of \cit!{2}, 3.8, by theorem~\nmb!{1} above.
\qed\enddemo

\enddocument